\documentclass{amsart}
\usepackage{amsmath}
\usepackage{amssymb}
\usepackage{amscd}
\usepackage{graphicx}
\usepackage{latexsym}
\usepackage{psfrag}
\usepackage{framed, color}

\usepackage{mathtools}

\newtheorem{thm}{Theorem}[section]
\newtheorem{lem}[thm]{Lemma}
\newtheorem{cor}[thm]{Corollary}

\theoremstyle{definition}
\newtheorem{defn}[thm]{Definition}

\theoremstyle{remark}
\newtheorem{rk}[thm]{Remark}

\makeatletter
 
 \@addtoreset{equation}{section}
\makeatother

\newcommand{\Int}{\mathop{\mathrm{Int}}\nolimits}
\newcommand{\rank}{\mathop{\mathrm{rank}}\nolimits}

\newcommand{\R}{\mathbb{R}}
\newcommand{\C}{\mathbb{C}}

\def\spmapright#1{\smash{%
 \mathop{\hbox to 1.3cm{\rightarrowfill}}
  \limits^{#1}}}
\def\spmapleft#1{\smash{%
 \mathop{\hbox to 1.3cm{\leftarrowfill}}
  \limits^{#1}}}

\begin{document}

\title[Relative trisection maps]{Relative trisection maps on $4$--manifolds}
\author[O.~Saeki]{Osamu Saeki}%
\address{Institute of Mathematics for Industry,
Kyushu University, Motooka 744, Nishi-ku, Fukuoka 819-0395,
Japan}

\date{\today}

\begin{abstract} 
We show that for an arbitrary compact orientable smooth $4$--dimensional
manifold $M$ with a given open book structure on the boundary $\partial M
\neq \emptyset$, which may be disconnected,
there exists a generic smooth map of $M$ onto the $2$--disk
which gives a relative trisection on $M$
consistent with the given open book structure on $\partial M$.
This gives a new proof to the existence of a relative trisection on $M$.
In fact, it turns out that every relative trisection can
be obtained in this way.
As a byproduct, we also show that $M$ always admits a decomposition
into a PALF (positive allowable Lefschetz fibration) over the $2$--disk
$D^2$ and a fiber
bundle over the circle with a sutured compression body fiber, consistent
with the given open book structure on $\partial M$ and with
the surface bundle structure of the PALF over $\partial D^2$.
\end{abstract}

\maketitle

\pagestyle{headings} 

\section{Introduction}\label{section1}

In \cite{GK4}, Gay and Kirby introduced the notion
of a trisection of a smooth closed connected orientable $4$--dimensional
manifold and showed that every such $4$--manifold admits a
trisection. A \emph{trisection} is a decomposition of a given
$4$--manifold into three $4$--dimensional handlebodies
which are attached along their boundaries in such a way
that the attaching parts correspond to standard
Heegaard splittings of the handlebody boundaries.
For the existence of a trisection, Gay--Kirby gave
two proofs: one is based on Morse $2$--functions
(or generic maps)
of $4$--manifolds into $\R^2$ \cite{GKPNAS, GK2}, and the other one
is based on the usual handle decomposition of $4$--manifolds.
Gay and Kirby also considered trisections of
$4$--manifolds with boundary, where they implicitly
assumed that the boundary was connected.
Later, the notion of trisection for closed $4$--manifolds
was extended to $4$--manifolds with (possibly disconnected) boundary by
Castro, Gay and Pinz\'{o}n-Caicedo \cite{CGP}, and
it is now called a \emph{relative trisection}.
An interesting feature of a relative trisection
is that such a structure naturally induces an
open book structure on the boundary $3$--manifold,
and they showed that for a given smooth compact
connected orientable $4$--manifold such that
an open book structure is given on the boundary,
there always exists a relative trisection
consistent with the given open book structure on the boundary.
Their proof was essentially based on a handlebody
argument and they did not give a proof based on
Morse $2$--functions into $\R^2$.

The main purpose of this paper is to show that
for a given smooth compact
connected orientable $4$--manifold possibly with
disconnected boundary such that
an open book structure is given on the boundary,
there always exists a Morse $2$--function on
the $4$--manifold onto the $2$--disk which
is compatible with the open book structure
given on the boundary and which gives a relative
trisection. Here, the open book structure on each
component of the boundary $3$--manifold may possibly have
an empty binding.
Such a Morse $2$--function
is called a \emph{relative trisection map}.
For a precise definition, the reader is referred to
Section~\ref{section2} of the present paper.

Our strategy for the proof is to extensively use 
the realizable moves for Morse $2$--functions (or
generic maps) with Lefschetz singularities
developed by Baykur and the author \cite{BSPNAS, BS1},
and later re-organized in \cite{Sae25} based on
Cerf theory \cite{Cerf, HW}.
In fact, it has been proved in \cite{BS1} that
every smooth closed connected orientable $4$--manifold
admits a \emph{simplified} trisection map, where
a trisection map is simplified if it is injective
when restricted to the singular point set.
In our case, because of the presence of nonempty boundary,
the same method does not directly apply and
we need certain new ideas. This will be thoroughly
explained in Section~\ref{section3}.

In Section~\ref{section3.5}, we show that
a given relative trisection of a $4$--manifold is
always realized by a relative trisection map
as long as the binding of the open book structure
on each boundary component is nonempty: such
a decomposition is obtained by
pulling back the three sectors of $D^2$ by
an appropriate relative trisection map onto $D^2$.

Unfortunately, our relative trisection map
will not be simplified in general. However, in the course
of the construction in Section~\ref{section3}, we get a nice map
onto the $2$--disk which has an interesting
property. In Section~\ref{section4}, we will see that
this leads to a decomposition
of a smooth compact connected oriented $4$--manifold
such that an open book structure with nonempty binding is given on the boundary,
into two pieces; one is a sutured compression body
bundle over $S^1$ extending the open book structure
on the boundary and the other is a PALF
(positive allowable Lefschetz fibration).
As corollaries we also show that such a $4$--manifold
always admits a broken Lefschetz fibration \cite{ADK, BSPNAS, BS1, GK1}
consistent with the given open book structure
on the boundary, and that any orientation preserving
diffeomorphism of a compact oriented surface
is ``bordant'' to a product of positive Dehn twists
through a self-diffeomorphism of a sutured compression body.

Throughout the paper, all manifolds and maps are
smooth of class $C^\infty$ unless otherwise stated.
The symbol ``$\cong$'' denotes a diffeomorphism
between manifolds. 
For a smooth map $f : M \to M'$ between smooth manifolds,
$S(f) (\subset M)$ denotes the set of singular points of $f$.

\section{Preliminaries}\label{section2}

In this section, we recall some necessary materials
which will be used in the following sections.

\begin{defn}\label{def:openbook1}
Let $N$ be a closed $n$--dimensional manifold.
Suppose that $B$ is a closed submanifold of $N$ of dimension $n-2$ 
with trivial normal bundle, and let $\tau : T(B) \to B \times D^2$
be a trivialization of the normal disk bundle
$T(B)$ of $B$ in $N$, where $T(B)$ is identified with a closed
tubular neighborhood of $B$ in $N$. 
We assume that there exists a smooth
$S^1$--bundle $\varphi : N \setminus B \to S^1$
such that $\varphi|_{T(B) \setminus B}$
coincides with the composition 
\begin{equation}
T(B) \setminus B \spmapright{\tau|_{T(B) \setminus B}} B \times (D^2
\setminus \{0\}) \spmapright{pr_2} D^2 \setminus \{0\} 
\spmapright{r} S^1,
\label{eq:book}
\end{equation}
where $pr_2$ is the projection to the second factor and
$r$ is the radial projection.
Such a structure $OB = (N, B, T(B), \varphi)$ of $N$ is called 
an \emph{open book structure}.
Note that for each $\theta \in S^1$,
the closure $P_\theta$ of $\varphi^{-1}(\theta)$ in $N$ coincides
with $\varphi^{-1}(\theta) \cup B$ and is a compact
$(n-1)$--dimensional submanifold of $N$ with boundary $B$.
The submanifold $B$ is called the \emph{binding}
and each $P_\theta$ (or $P_\theta \setminus \Int{T(B)}$)
is called a \emph{page} of the open book
structure. Note that we allow $B$ to be empty: in such
a case, $N$ is the total space of a smooth fibration
over the circle with fiber a closed $(n-1)$--dimensional
manifold.
\end{defn}

For some historical account of open book
structures, the reader is referred to
Quinn \cite[\S2]{Quinn} and Myers \cite{My}.
In particular, it has been known that every closed orientable
$3$--dimensional manifold admits an open book structure
\cite{Alex}.

Suppose now that $M$ is a compact connected
orientable $4$--dimensional manifold such that
an open book structure $OB = (N, B, T(B), \varphi)$ is given
on $N = \partial M$. Note that the binding
$B$ is a link in $\partial M$.
Let $f: M \to D^2 (\subset \R^2 = \C)$ be a smooth map and
set $g = f|_{\partial M} : \partial M \to D^2$.
We assume that $g$ is \emph{associated with $OB$} in the
following sense.
The set $g^{-1}(\partial D^2) = N \setminus \Int{T(B)}$ is a compact $3$--dimensional
submanifold of $\partial M$ on which $g$ coincides with the fibration
$\varphi|_{N \setminus \Int{T(B)}} : N \setminus \Int{T(B)} \to
S^1 = \partial D^2$, and $T(B)$ is the closure
of $g^{-1}(\Int{D^2})$ and is a disjoint union of solid tori 
$S^1 \times D^2$ on each of which 
$g$ is the projection $S^1 \times D^2 \to D^2$.
The binding $B$ of $OB$ is $g^{-1}(0)$ and the page over 
$\theta \in S^1 = \partial D^2$ is $P_\theta = 
g^{-1}(\{\lambda \theta \,|\,
0 \leq \lambda \leq 1\})$.
The monodromy is the isotopy class (rel.\ boundary)
of the return map $h : P_1 \to P_1$,
$1 \in S^1 \subset \C$,
associated with
any vector field transverse to the interiors of all the
pages and meridional near the binding.

Furthermore, we assume that $f^{-1}(\partial D^2) \subset \partial M$ 
and $f$ is a submersion on a 
neighborhood of $\partial M$.
In this case, we say that $f$ is \emph{compatible} with the
open book structure $OB$.

Let $p \in \Int{M}$ be a \emph{singular point} of $f$, i.e.\ 
$\rank{df_p} < 2$. 
The map $f$ is said to have a \emph{fold singularity} at $p$ \cite{E1, E2}
if there are local coordinates 
around $p$ and $f(p)$ in which the map is given by
$$(t, x_1, x_2, x_3) \mapsto (t, \pm x_1^2 \pm x_2^2 \pm x_3^2),$$
and a \emph{cusp singularity} if the map is locally given by
$$(t, x_1, x_2, x_3) \mapsto (t, x_1^3 + t x_1 \pm x_2^2 \pm x_3^2).$$
Note that fold singularities constitute a $1$--dimensional closed submanifold
of $\Int{M}$ and that the cusp singularities constitute a finite set of points.
Furthermore, each cusp singularity is adjacent to the fold arc.
A fold or a cusp point $p$ is \emph{definite} if the coefficients of all quadratic terms 
in the corresponding local model as above are of the same sign, \emph{indefinite} otherwise. 
Note that an indefinite cusp is always adjacent to indefinite fold arcs
and that a definite cusp is adjacent to both definite and indefinite fold
arcs.

We see easily that the map restricted to fold arcs are immersions.
The immersed image of any indefinite (or definite) fold arc or circle is normally 
oriented in such a way that the normal orientation
indicates the direction in which the topology of a fiber changes by a $2$--handle
(resp.\ $3$--handle) attachment when crossing over the fold from one side to the other. 
This means that, for indefinite fold image,
a generic fiber over the region the arrow starts, if connected, has one higher genus
than the generic fiber over the region the arrow points into.

A special case of Thom's transversality implies that any smooth map from 
an $n \geq 2$ dimensional manifold to a surface can be approximated arbitrarily 
well by a map with only fold and cusp singularities \cite{Lev0, T, Wh}. 
If the above smooth map $f : M \to D^2$ 
has only fold and cusp singularities, then
it is called a \emph{generic map} (or an \emph{excellent map},
or --more recently-- a \emph{Morse $2$--function}). The \emph{singular point set} $S(f)$ 
of $f$ is contained in $\Int{M}$ and it is a disjoint union of finitely many circles, which are 
composed of finitely many cusp points, and arcs and circles of fold singularities. 
We call $f$ an \emph{indefinite generic map} if all of its fold singularities are indefinite.

Now, suppose that the $4$--dimensional manifold $M$ is oriented.
The map $f$ is said to have a \emph{(positive)
Lefschetz singularity} at a point 
$p \in \Int{M}$ if there are orientation preserving local coordinates 
(in the $C^\infty$ sense) around $p$ and $f(p)$ 
so that $f$ conforms to the complex local model
$$(z_1, z_2) \mapsto z_1 \, z_2.$$
If it is locally given by
$$(z_1, z_2) \mapsto z_1 \, \overline{z_2},$$
then it is called a \emph{negative (or achiral)
Lefschetz singularity}.
A \emph{broken Lefschetz fibration} (\emph{BLF}, for short)
is a surjective map $f: M \to D^2$ for which $S(f)$
consists of finitely many (positive) Lefschetz singularities 
and indefinite fold circles. 
In particular, a broken Lefschetz fibration does
not have any cusp singularities.
An \emph{achiral broken Lefschetz fibration}
(\emph{ABLF}, for short)
is a surjective map $f: M \to D^2$ for which $S(f)$ consists
of finitely many positive and negative
Lefschetz singularities and indefinite fold circles. 

Since the set of generic maps is open and dense in an appropriate mapping space 
endowed with the Whitney $C^{\infty}$ topology, every achiral
broken Lefschetz fibration can 
be approximated (and hence homotoped) to a map with only fold and cusp singularities. 
As shown in \cite{Lek} (see also \cite{BS1}), one can trade an indefinite cusp 
point with a (positive or negative)
Lefschetz singularity, and locally perturb a (positive or negative)
Lefschetz singularity into a simple indefinite 
singular circle with three cusp points, allowing one to switch between 
indefinite generic maps and 
achiral broken Lefschetz fibrations in a rather standard way.

With these in mind, we call a smooth surjective map $f: M \to D^2$ an 
\emph{indefinite fibration} if it is an indefinite generic map 
outside of a finite collection of (positive or negative)
Lefschetz singular points. 
For an indefinite fibration $f$, we call the 
closure $Z_f$ of the $1$--dimensional fold arcs 
and circles, containing cusp points,
the \emph{round locus}, 
and its image $f(Z_f)$ the \emph{round image} of $f$. 
Note that the restriction of an indefinite fibration to its 
round locus is an immersion except at the cusp points. 
Its round image is, generically, a collection
of cusped immersed curves in 
$\Int{D^2}$ with transverse double points off the cusp points, which we will 
assume to be the case hereon. 
We can moreover assume that there is at most one (positive or
negative) Lefschetz singular point 
on any fiber, and 
also, (positive or negative)
Lefschetz singular values and round image are disjoint.

Let $f : M \to D^2$ be a smooth map 
compatible with
an open book structure $OB$ given on $\partial M$.
We say that $f$ is a \emph{Lefschetz fibration}
(\emph{LF}, for short) if $f$ has at most positive
Lefschetz singularities, and $f$ is an
\emph{achiral Lefschetz fibration} (\emph{ALF},
for short) if $f$ has at most positive and/or negative
Lefschetz singularities.

Let $f : M \to D^2$ be an indefinite fibration.
We say that $f$ (or the round image $f(Z_f)$) is 
\emph{outward-directed},
if the complement of a regular value
$z_0 \in \Int{D^2}$ can be non-singularly foliated by arcs oriented from 
$z_0$ to $\partial D^2$, which
intersect the image of each fold arc (excluding the image of the
cusp points) transversely in its normal direction.

\begin{figure}[t]
\centering
\includegraphics[width=0.95\linewidth,height=0.3\textheight,
keepaspectratio]{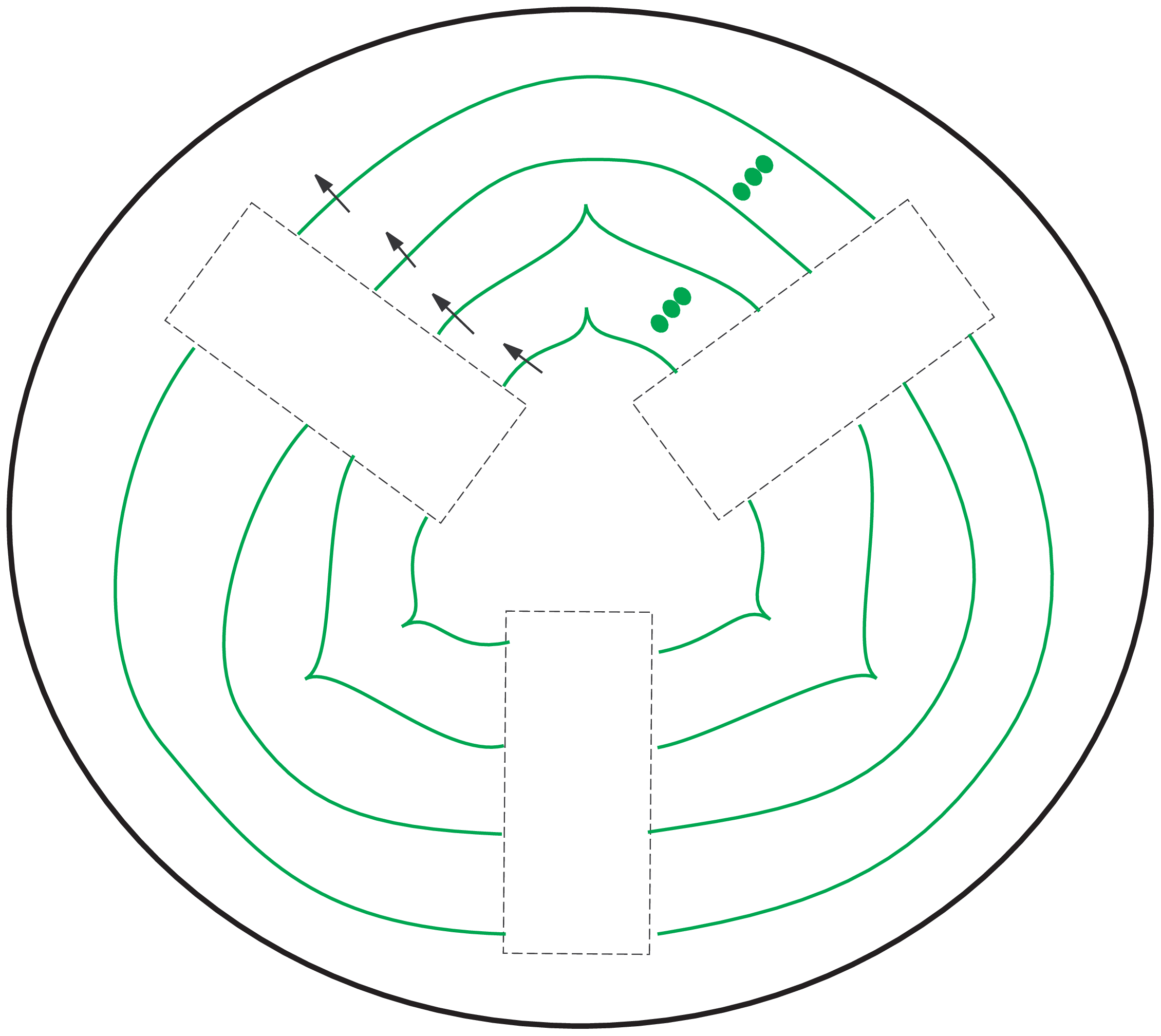}
\caption{Singular value set of a relative trisection map. The three boxes
represent Cerf graphics.
In this figure and the following, green curves depict the image
of indefinite folds and indefinite cusps, and
small arrows depict the relevant normal directions.}
\label{fig12}
\end{figure}

Let $f : M \to D^2$ be an indefinite generic map compatible with
an open book structure $OB$ on $\partial M$.
We say that $f$ is a \emph{relative trisection map}
if $f(S(f))$ is outward directed in $D^2$
and is situated as depicted in Fig.~\ref{fig12}, where
each box represents a Cerf graphic without cusps in the sense of \cite{GK4}
and the green lines represent the image of the singular point set.

The terminology ``relative trisection map''
seems to have first appeared in \cite{Olsen}, although
it is assumed that the boundary is connected.

\section{Construction of a relative trisection map}\label{section3}

Our main theorem of this paper is the following.

\begin{thm}\label{thm1}
Let $M$ be a compact connected orientable $4$--dimensional
manifold with nonempty boundary $\partial M$, which
may be disconnected.
Suppose an open book structure $OB$ on $\partial M$
is given. Then, there exists a relative
trisection map $f : M \to D^2$ compatible with $OB$.
\end{thm}

\begin{figure}[t]
\centering
\includegraphics[width=0.95\linewidth,height=0.3\textheight,
keepaspectratio]{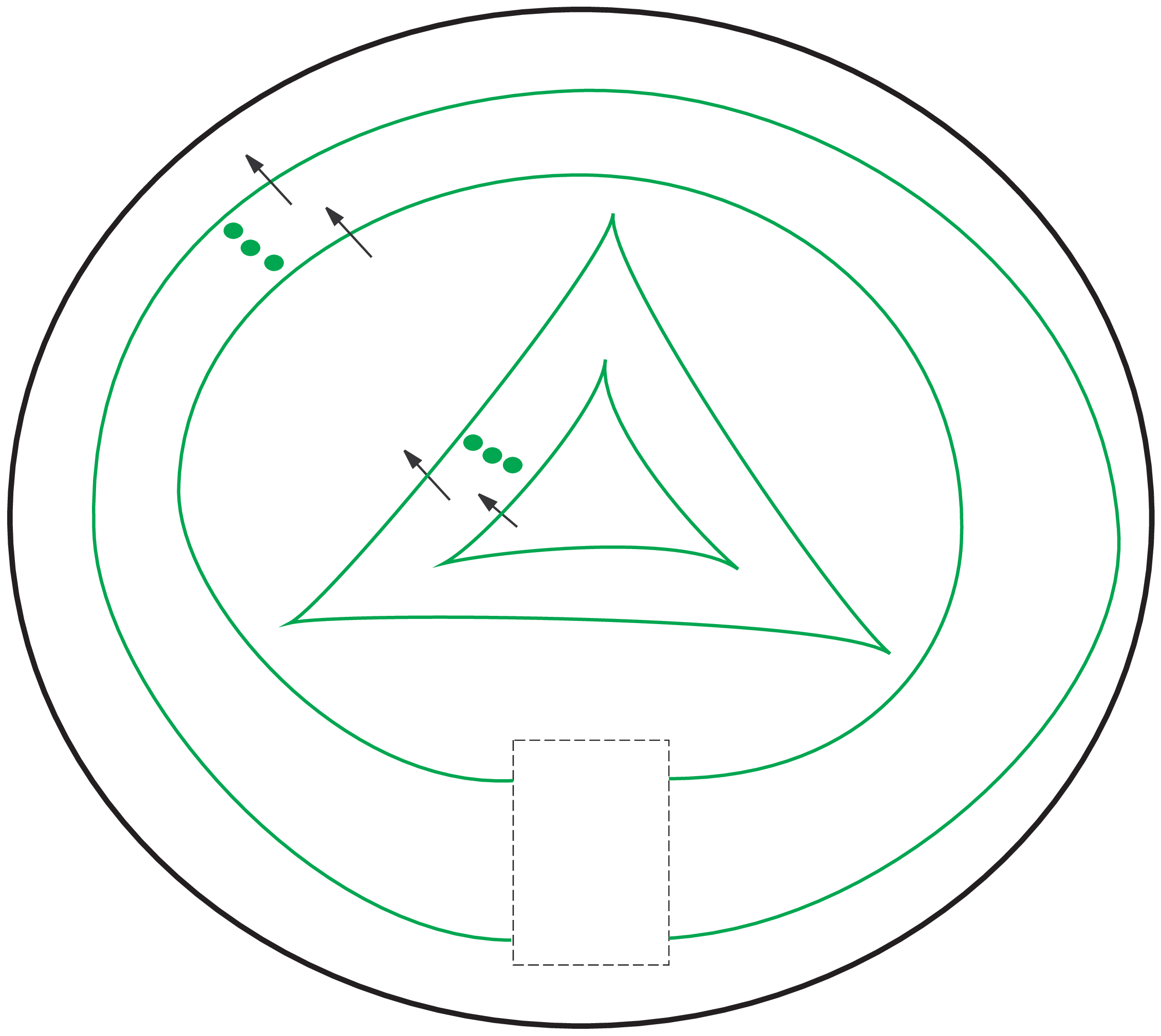}
\caption{Singular value set of a relative trisection map
constructed in the proof of Theorem~\ref{thm1}.}
\label{fig13}
\end{figure}

In fact, we will show that there exists a
relative trisection map whose singular value
set is as depicted in Fig.~\ref{fig13} with
a box representing a Cerf graphic without cusps, but
possibly with crossings.

\begin{proof}[Proof of Theorem~\textup{\ref{thm1}}]
By the given open book structure $OB = (\partial M,
B, T(B), \varphi)$ on $\partial M$, we have the natural decomposition
$\partial M = T(B) \cup Q$, where $Q = \partial M \setminus \Int{T(B)}$.
We may assume that
$M$ has corners exactly along $\partial T(B) = \partial Q$.
Note that $T(B)$ is a trivial disk bundle neighborhood of the binding $B$
and that $Q$ fibers over the circle $S^1 = \partial D^2$.
Let $f_1 : \partial M \to D^2$ be a smooth map such that
$(f_1)^{-1}(\partial D^2) = Q$, that $f_1|_Q : Q \to \partial D^2$
coincides with $\varphi|_Q : Q \to S^1$, and that
$f_1|_{T(B)}: T(B) \to D^2$ is the composition $pr_2 \circ \tau$,
where $\tau : T(B) \to B \times D^2$ is the trivialization
and $pr_2 : B \times D^2 \to D^2$ is the projection to the
second factor.
Then, we can extend $f_1$ to a generic map
$M \to D^2$, still denoted by $f_1$, which is compatible with $OB$,
since $D^2$ is contractible. Note that the singularities
of $f_1$ are contained in $\Int{M}$ and are,
in general, composed of definite folds, indefinite
folds, definite cusps and indefinite cusps.

Let $f_2 : M \to S^2$ be a smooth map which
is the composition of $f_1 : M \to D^2$ and
an embedding $D^2 \to S^2$.
Then, by using techniques of \cite{Sa, Sa2}, we can
homotopically modify $f_2$ on $\Int{M}$
to get a map $f_3: M \to S^2$ such that
$f_3|_{\partial M} = f_2|_{\partial M} : \partial M
\to D^2 \subset S^2$ and $f_3|_{\Int M} : \Int{M}
\to S^2$ is a generic map without definite fold
singularities. Note that in our case, we have
nonempty boundary for $M$; however, we see easily
that the same procedure as in \cite{Sa2}
works well. So, $f_3$ is a map with only
indefinite folds and indefinite cusps as its
singularities.

Now, by modifying $f_3$ by homotopy on a collar
neighborhood of $\partial M$ in such a way that
the image of $\partial M$ engulfs
the singular point set image, we get a map
$f_4 : M \to S^2$ such that $\Delta_1 =
f_4(\partial M) \cong D^2$ and
$f_4(S(f_4)) \subset \Int \Delta_1$
(see Fig.~\ref{fig1}).
Note that $f_4|_{\partial M} : \partial M
\to \Delta_1 \cong D^2$ can be identified with
$f_1|_{\partial M} : \partial M \to D^2$.
Set $\Delta_2 = S^2 \setminus \Int \Delta_1
\cong D^2$ and $V = (f_4)^{-1}(\Delta_2) \setminus \partial M$.
Note that $f_4|_V : V \to \Delta_2$ is a trivial fiber bundle with fibers being
closed surfaces, which may be disconnected.

\begin{figure}[t]
\centering
\psfrag{f3}{$f_3(S(f_3))$}
\psfrag{f4}{$f_4(S(f_4))$}
\psfrag{D}{$D^2 = f_3(\partial M)$}
\psfrag{D1}{$\Delta_1 = f_4(\partial M) \cong D^2$}
\psfrag{D2}{$\Delta_2$}
\includegraphics[width=0.95\linewidth,height=0.4\textheight,
keepaspectratio]{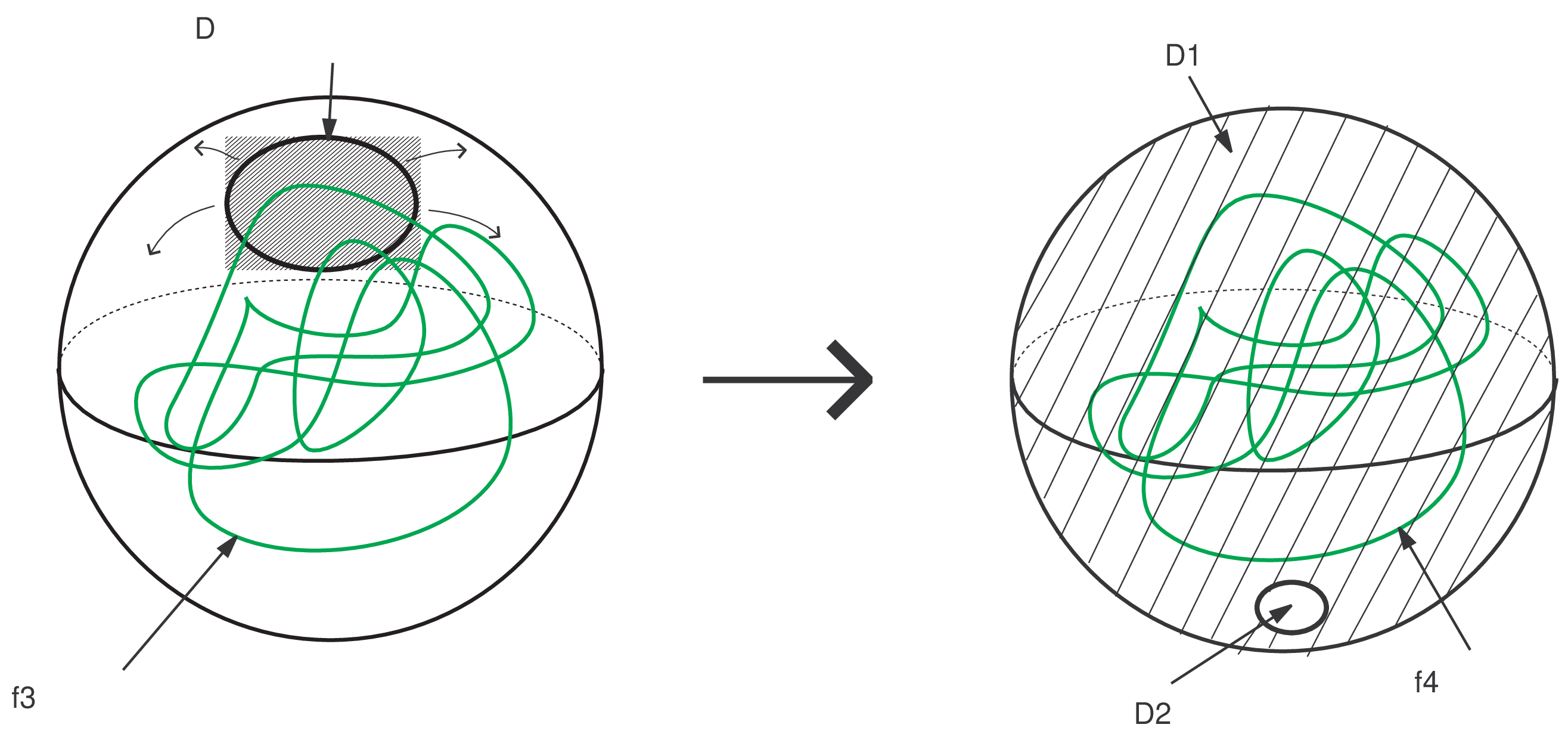}
\caption{Modifying $f_3$ on a neighborhood of $\partial M$.}
\label{fig1}
\end{figure}

We can now modify $f_4$ on $(f_4)^{-1}(\Int \Delta_1) \cap \Int M$
to get a map $f_5 : M \to S^2$
which has only indefinite folds and Lefschetz singularities
such that $f_5(S(f_5)) \subset \Int \Delta_1$ is outward
directed (see \cite[Theorem~4.1 and Remark~4.5]{BS1}).
Refer to Fig.~\ref{fig2}.

\begin{figure}[t]
\centering
\psfrag{f5p}{$f_5(\partial M) = \Delta_1 \subset S^2$}
\psfrag{f5}{$f_5(S(f_5))$}
\includegraphics[width=0.95\linewidth,height=0.3\textheight,
keepaspectratio]{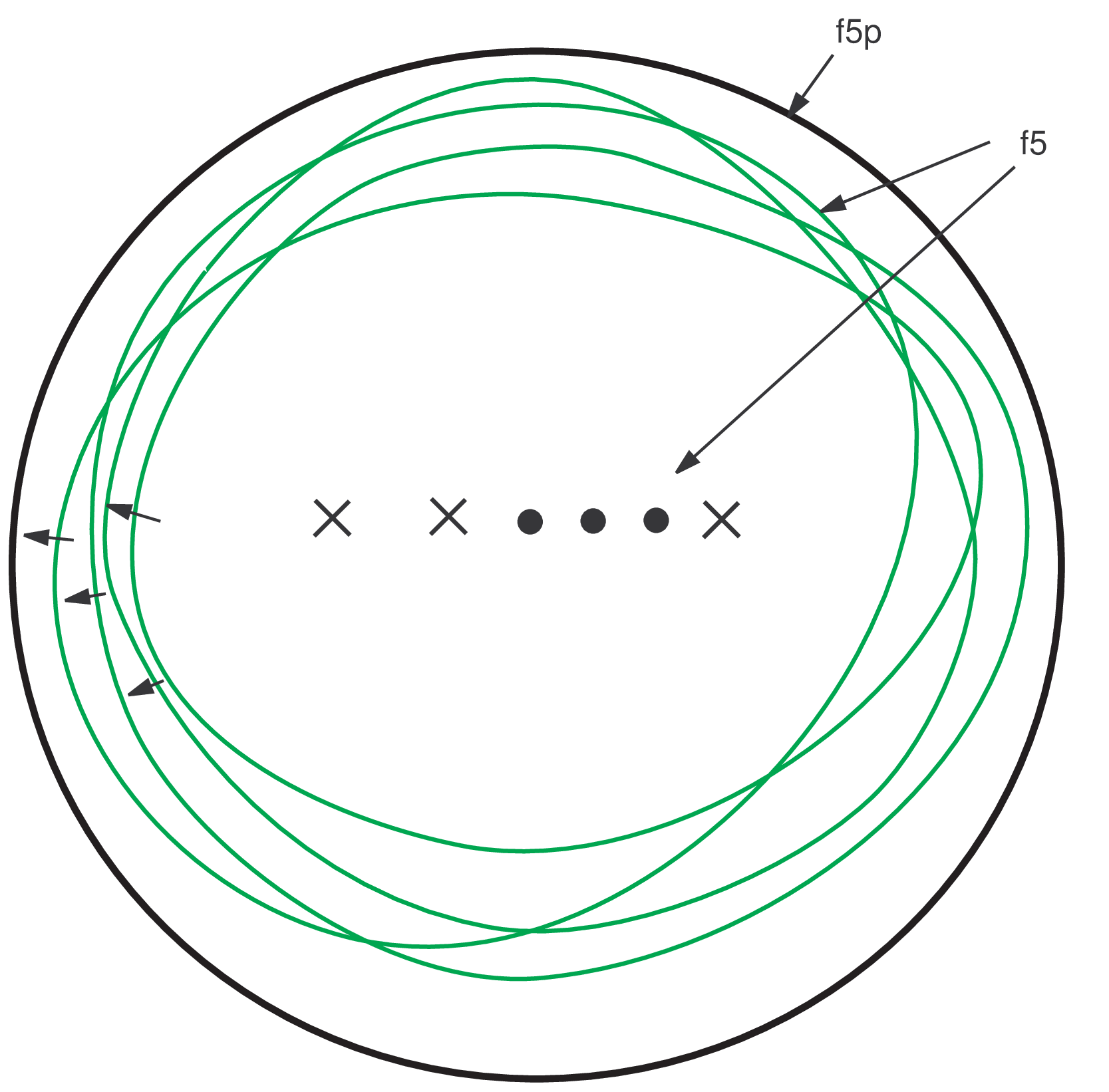}
\caption{The singular point set image of $f_5 : M \to S^2$. In this
figure and the following, 
small crosses ``$\times$'' depict the images of Lefschetz singularities.}
\label{fig2}
\end{figure}

By modifying $f_5$ on a neighborhood of $V$ by
the same procedure as in \cite[Theorem~7.1]{BS1},
we get a map
$f_6 : M \to \Delta_1 \cong D^2$ which
has definite folds, indefinite folds and Lefschetz
singularities
such that
\begin{enumerate}
\item $f_6(\partial M) = f_6(M) = \Delta_1$,
\item $f_6|_{\partial M} : \partial M \to \Delta_1$
can be identified with $f_1|_{\partial M} : \partial M \to D^2$,
\item the image of the definite folds of $f_6$ consists
of a finite family of disjoint concentric circles
near $\partial \Delta_1$,
\item the image of the indefinite folds is contained
in the inner region, is outward directed and is
disjoint from the image of the definite folds,
\item the images of the Lefschetz singularities
are distinct, located in the innermost
region $\Delta_3$ in $\Int{\Delta_1}$ disjoint
from the definite and indefinite fold
images.
\end{enumerate}
Refer to Fig.~\ref{fig3}.
Note that in \cite[Theorem~7.1]{BS1}, the fibers in $V$
over $\Delta_2$ are connected, while in our case,
they may be disconnected. Therefore, we may have several
definite fold circles, one for each
component of a fiber.

\begin{figure}[t]
\centering
\psfrag{f6p}{$f_6(\partial M) = f_6(M) = \Delta_1$}
\psfrag{f6}{$f_6(S(f_6))$}
\psfrag{D3}{$\Delta_3$}
\psfrag{l}{$\ell$}
\psfrag{a}{$a$}
\psfrag{b}{$b$}
\includegraphics[width=0.95\linewidth,height=0.3\textheight,
keepaspectratio]{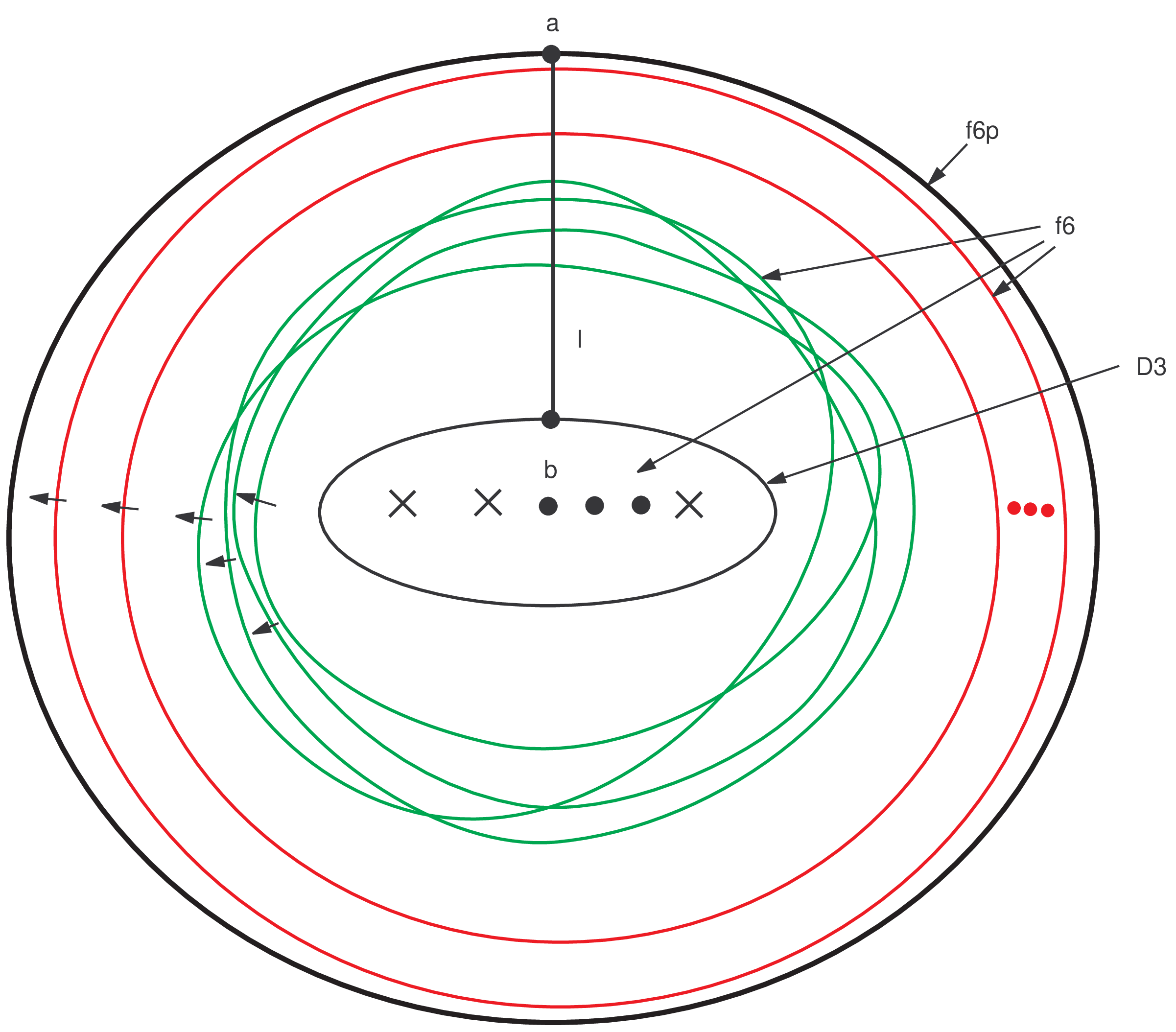}
\caption{The singular point set image of $f_6 : M \to \Delta_1$. In this
figure and the following, red curves depict the image of definite folds.}
\label{fig3}
\end{figure}

Since $\Delta_3$ deformation retracts
to a $1$--dimensional CW complex
consisting of a point $x_0 \in \Int{\Delta_3}$ distinct from the Lefschetz singularity
images and some embedded arcs which are disjoint except at $x_0$ and
which connect that point to the images of the
Lefschetz singularities, we see that $(f_6)^{-1}(\Delta_3)$
deformation retracts to a $2$--dimensional CW complex
in $M$. This implies that $M'$, the closure of $M \setminus
(f_6)^{-1}(\Delta_3)$, is connected, 
since so is $M$.
Note that $f_6|_{M'} : M' \to \Delta_1
\setminus \Int{\Delta_3}$ composed with a radial
projection $pr: \Delta_1
\setminus \Int{\Delta_3} \to \partial \Delta_1$
is a smooth fiber bundle. Let $\ell \subset 
\Delta_1 \setminus \Int{\Delta_3}$ be a fiber
of $pr$, which is a line segment, and set
$\partial \ell = \{a, b\}$, where $a \in \partial \Delta_1$
and $b \in \partial \Delta_3$.
The fiber $K$ of the fibration $pr \circ f_6|_{M'}: M'
\to \partial \Delta_1$ is a compact $3$--dimensional
manifold, and $f_6|_K : K \to \ell$
is a Morse function with $(f_6)^{-1}(a)$ a page
of the given $OB$. In other words, $K$ can
be considered to be a relative cobordism
between a page of $OB$ and the compact surface $(f_6)^{-1}(b)$.

If $K$ is disconnected, then we see that the
monodromy of the fibration over $\partial \Delta_1$
acts non-trivially on $\pi_0(K)$, which is a contradiction,
since the total space is connected and
the monodromy action on $\pi_0((f_6)^{-1}(a))$ is trivial
by the definition of an open book structure.
Consequently, $K$ is connected.

Now, let us consider the Morse function
$f_6|_K : K \to \ell$, where $\ell$
is identified with $[0, 1]$ in such a way
that $a$ (or $b$) is identified
with $0$ (resp.\ $1$). It has several
index $0$ critical points, each corresponding
to a definite fold component of $f_6$.
As $(f_6)^{-1}(a)$ is nonempty and $K$
is connected, we see that for an index $0$
critical point $p$, there is a canceling index $1$
critical point $q$, whose descending manifold $\delta$
``connects'' $p$ to a component of $(f_6)^{-1}(a)$.

Then, by the procedure as depicted in
\cite[Figure~19]{Sae25}, we can cut
the definite fold circle at $p$ with the
expense of creating two definite cusps
by applying the fold merge move. (See \cite[Figure~3]{BS1}.
The fold merge move is not ``always-realizable'' in
general; however, in our case of canceling pairs
of critical points, it is realizable.)
More precisely, we proceed as in Fig.~\ref{fig10},
using the moves $\mathrm{R2}^1$ of \cite{BS1}.

\begin{figure}[t]
\centering
\psfrag{p}{$f_6(p)$}
\psfrag{q}{$f_6(q)$}
\psfrag{D3}{$\Delta_1$}
\psfrag{d}{$f_6(\delta)$}
\psfrag{a}{$a$}
\includegraphics[width=0.95\linewidth,height=0.5\textheight,
keepaspectratio]{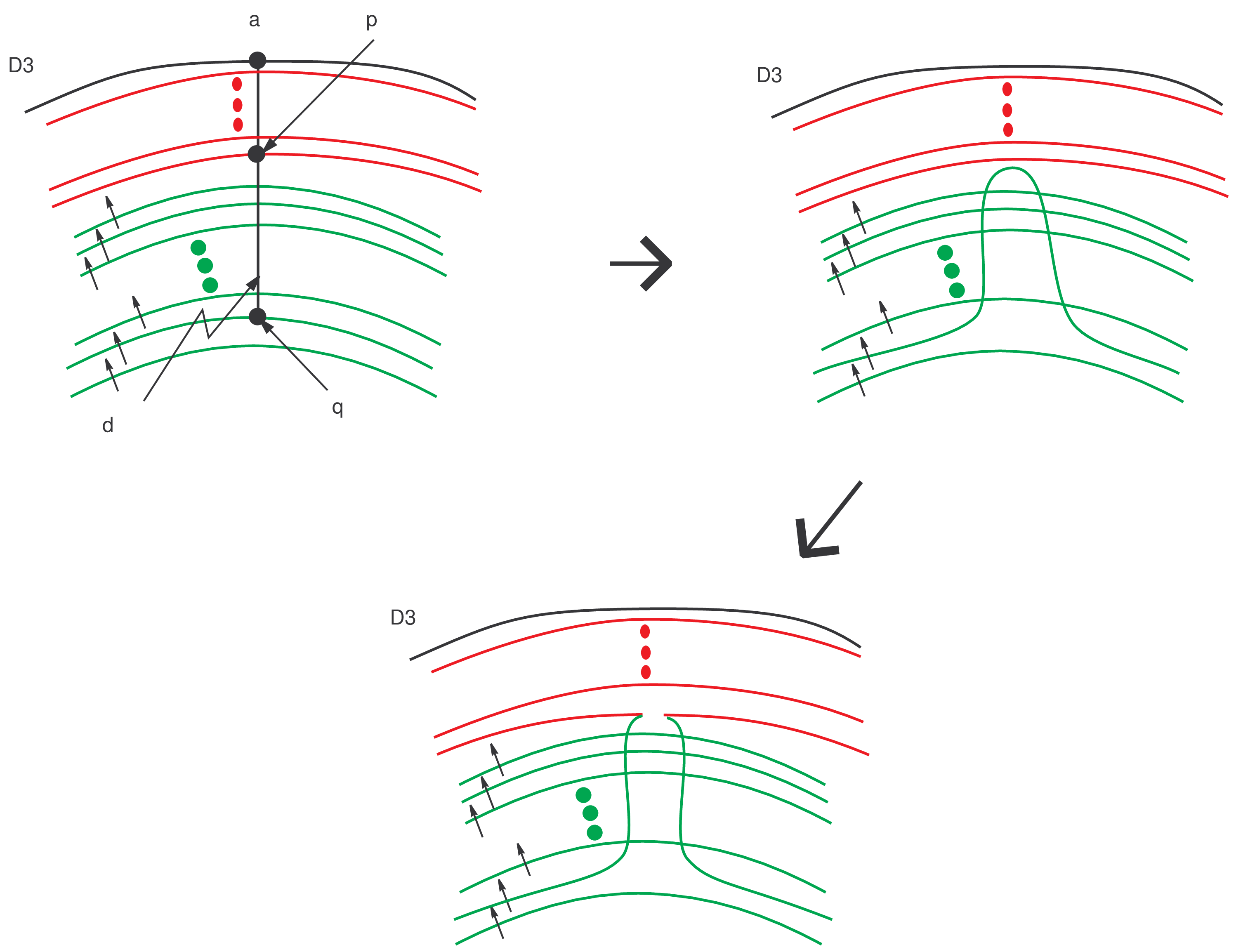}
\caption{We first move the image of the indefinite fold
arc near $q$ toward $f_6(p)$ and then we merge
the definite and indefinite fold arcs using $\delta$.}
\label{fig10}
\end{figure}

Iterating this procedure finitely many times,
we get a map $f_7 : M \to \Delta_1$
whose image of the definite fold singularities 
consists of disjoint embedded arcs in $\Int{\Delta_1}$
(see Fig.~\ref{fig4}).

\begin{figure}[t]
\centering
\psfrag{f7p}{$f_7(\partial M) = \Delta_1$}
\psfrag{f7}{$f_7(S(f_7))$}
\psfrag{D3}{$\Delta_3$}
\psfrag{l}{$\ell$}
\psfrag{a}{$a$}
\psfrag{b}{$b$}
\includegraphics[width=0.95\linewidth,height=0.3\textheight,
keepaspectratio]{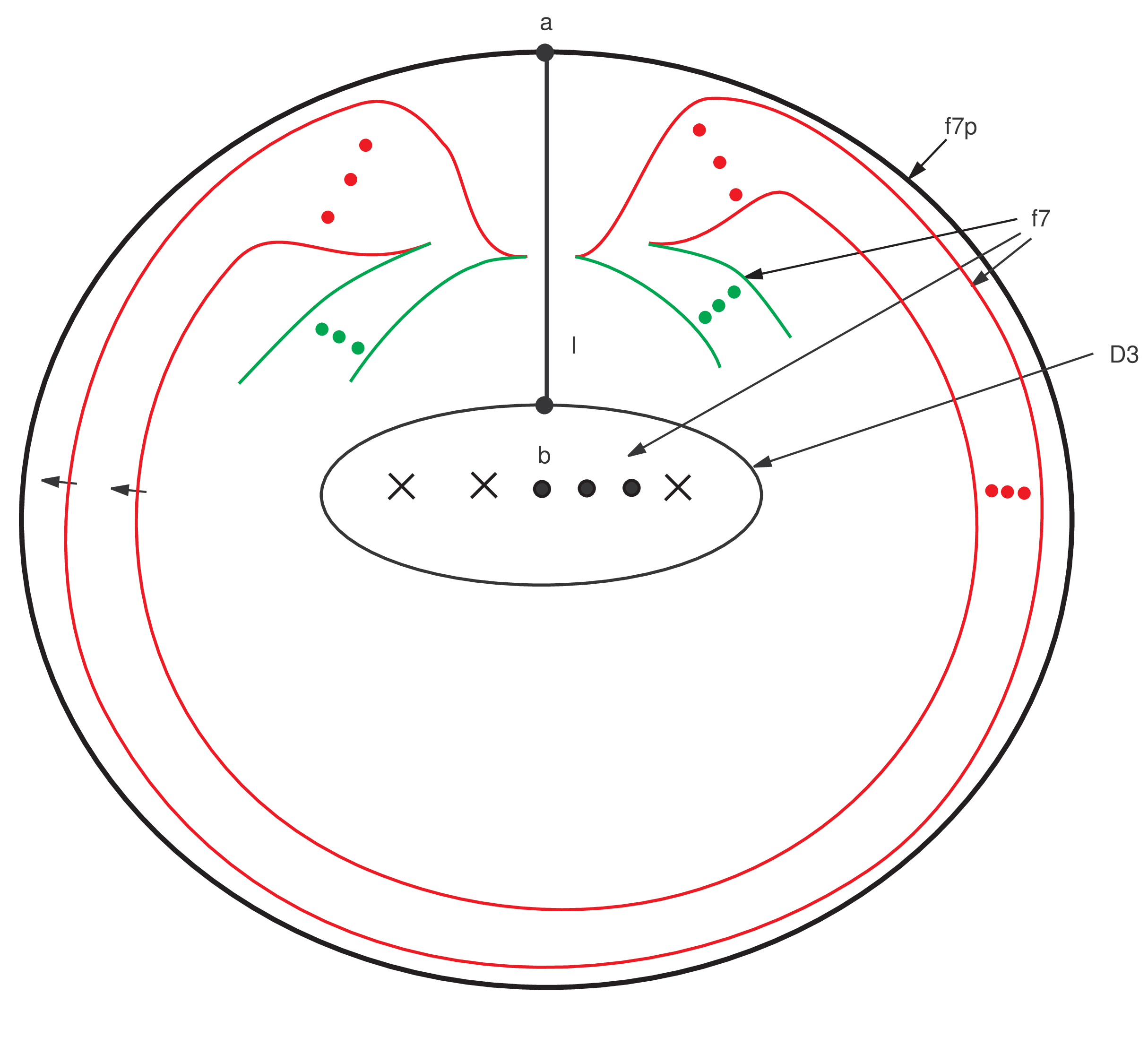}
\caption{By cutting each definite fold circle
by a fold merge move, we get $f_7 : M \to \Delta_1$
whose singular point set image is as depicted here.}
\label{fig4}
\end{figure}

Then, for each definite fold arc $\alpha$,
we can find a joining curve $\gamma$ for the pair of
definite cusp points adjacent to $\alpha$
in the sense of Levine \cite{Lev} such that
$f_7(\gamma)$ is ``parallel'' to $f_7(\alpha)$,
lying on the outward side.
Such a joining curve can be found, for the definite
fold arc $\alpha$ originating from the definite fold
circle containing the critical point $p$ of index $0$
as above, in such a way that it is in a small collar
neighborhood of the component of
$\partial M$ containing the component of
$(f_6)^{-1}(a)$ connected to $p$ by $\delta$.
Then, by the Levine cusp elimination technique
\cite{Lev} applied to $\gamma$,
we can eliminate the two definite cusps adjacent
to the definite fold arc $\alpha$ in such a way
that the definite fold arc turns into an ``inward'' circle
(see Fig.~\ref{fig5}). 

\begin{figure}[t]
\centering
\psfrag{f7p}{$f_7(\partial M) = \Delta_1$}
\psfrag{f7}{$f_7(S(f_7))$}
\psfrag{D3}{$\Delta_3$}
\psfrag{l}{$\ell$}
\psfrag{a}{$a$}
\psfrag{b}{$b$}
\psfrag{g}{$f_7(\gamma)$}
\psfrag{al}{$f_7(\alpha)$}
\includegraphics[width=0.95\linewidth,height=0.5\textheight,
keepaspectratio]{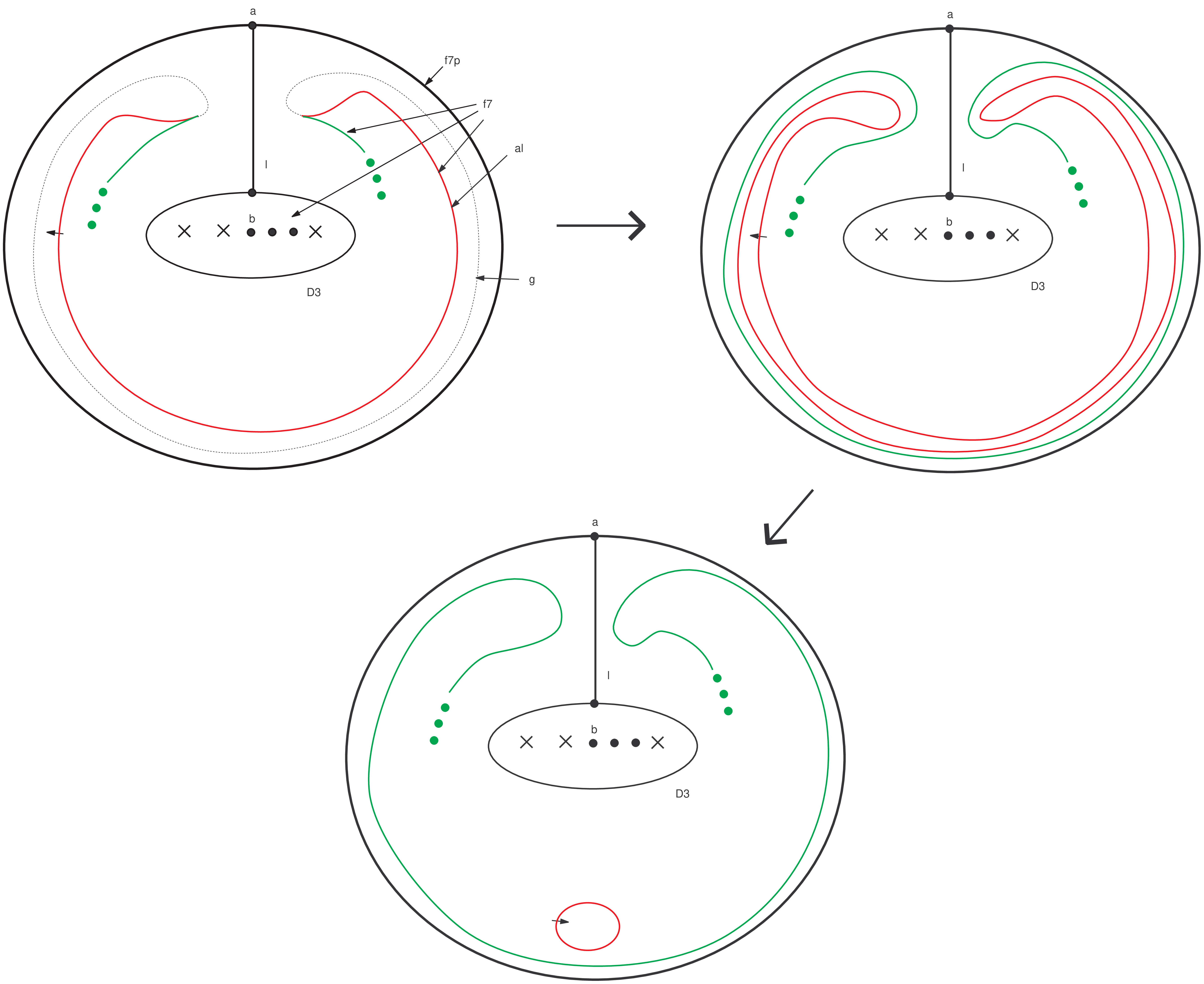}
\caption{Eliminating two definite cusps by pairs using
an appropriate joining curve $\gamma$, whose image is
depicted by a dotted curve.}
\label{fig5}
\end{figure}

Then, by the procedure as in \cite[Figure~7]{GKPNAS},
we can eliminate the definite fold circle
with the expense of creating an indefinite fold
circle in such a way that the image is embedded
and is outward directed.
We perform the same procedure for each
of the definite fold arc and get
an indefinite fibration $f_8 : M \to
\Delta_1$ (see Fig.~\ref{fig6}).

\begin{figure}[t]
\centering
\psfrag{f8p}{$f_8(\partial M) = \Delta_1$}
\psfrag{f8}{$f_8(S(f_8))$}
\psfrag{D3}{$\Delta_3$}
\includegraphics[width=0.95\linewidth,height=0.3\textheight,
keepaspectratio]{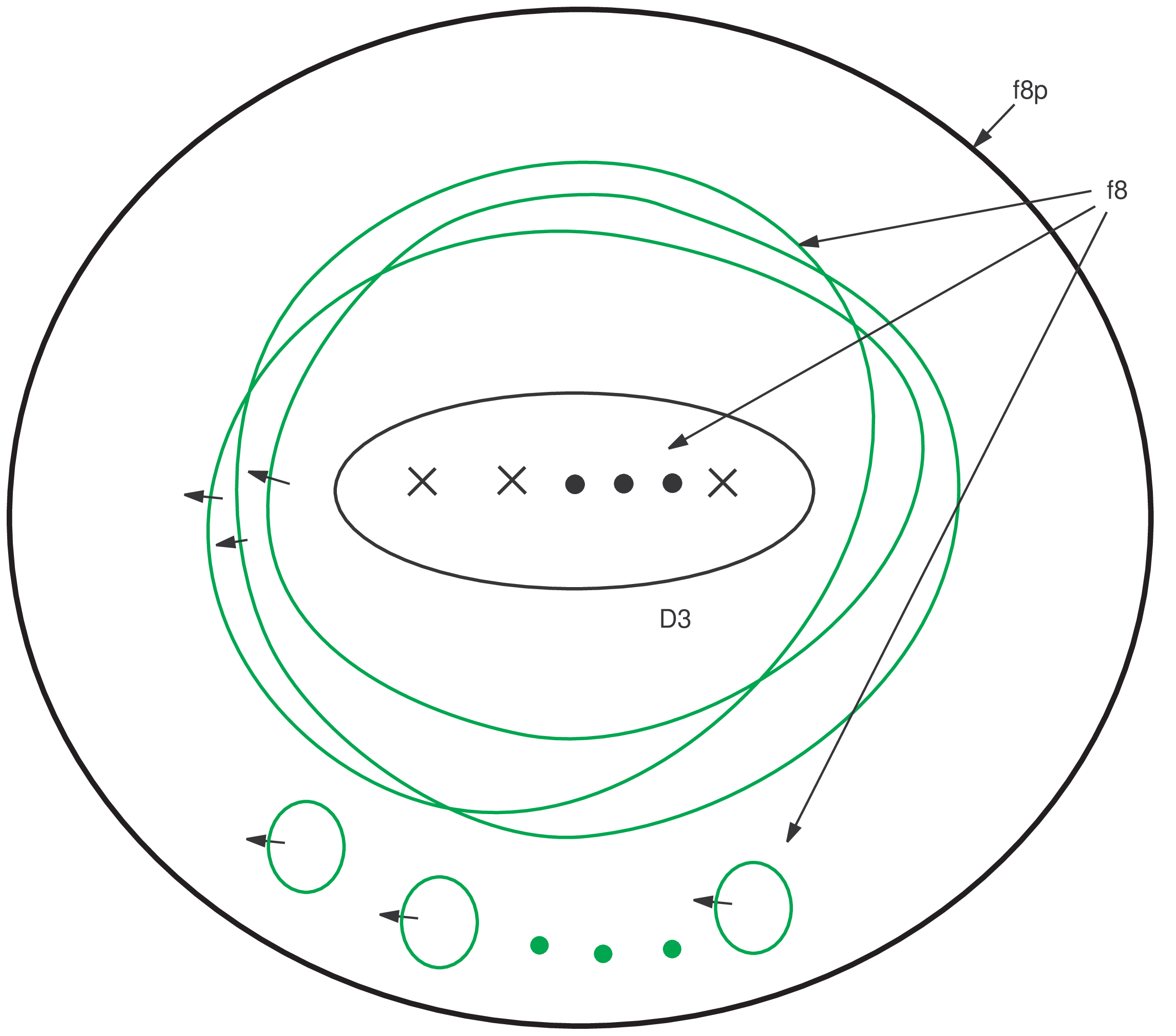}
\caption{Replacing each inward definite fold circle image
by an outward indefinite fold circle image.}
\label{fig6}
\end{figure}

Then, by applying $\mathrm{R2}^0$ moves together with
push moves as in \cite{BS1}, we can
move the image of the indefinite folds
so that it is outside of $\Delta_3$ and is outward directed,
and that the images of the Lefschetz singularities
are all distinct and are situated in $\Delta_3$.
Let this map
be denoted by $f_9 : M \to \Delta_1$
(see Fig.~\ref{fig7}).

\begin{figure}[t]
\centering
\psfrag{f9p}{$f_9(\partial M) = \Delta_1$}
\psfrag{f9}{$f_9(S(f_9))$}
\psfrag{D3}{$\Delta_3$}
\includegraphics[width=0.95\linewidth,height=0.3\textheight,
keepaspectratio]{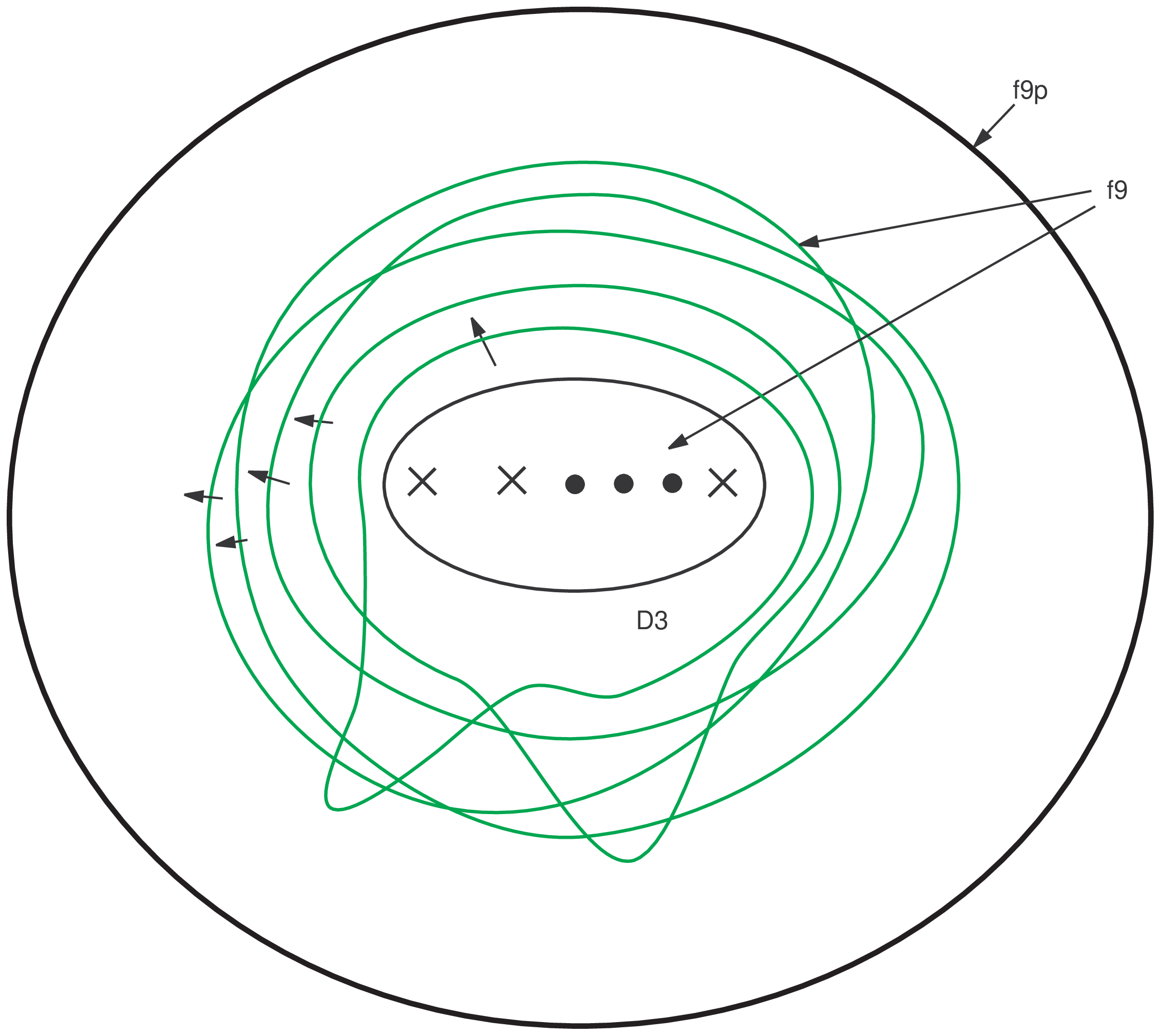}
\caption{We move each outward indefinite circle image
toward the center and then push the Lefschetz
critical point images inside. We do such
a process one by one for each outward
indefinite fold circle image as in Fig.~\ref{fig6}.}
\label{fig7}
\end{figure}

\begin{figure}[t]
\centering
\psfrag{f10p}{$f_{10}(\partial M) = \Delta_1$}
\psfrag{f10}{$f_{10}(S(f_{10}))$}
\psfrag{c}{$\cong$}
\includegraphics[width=0.95\linewidth,height=0.7\textheight,
keepaspectratio]{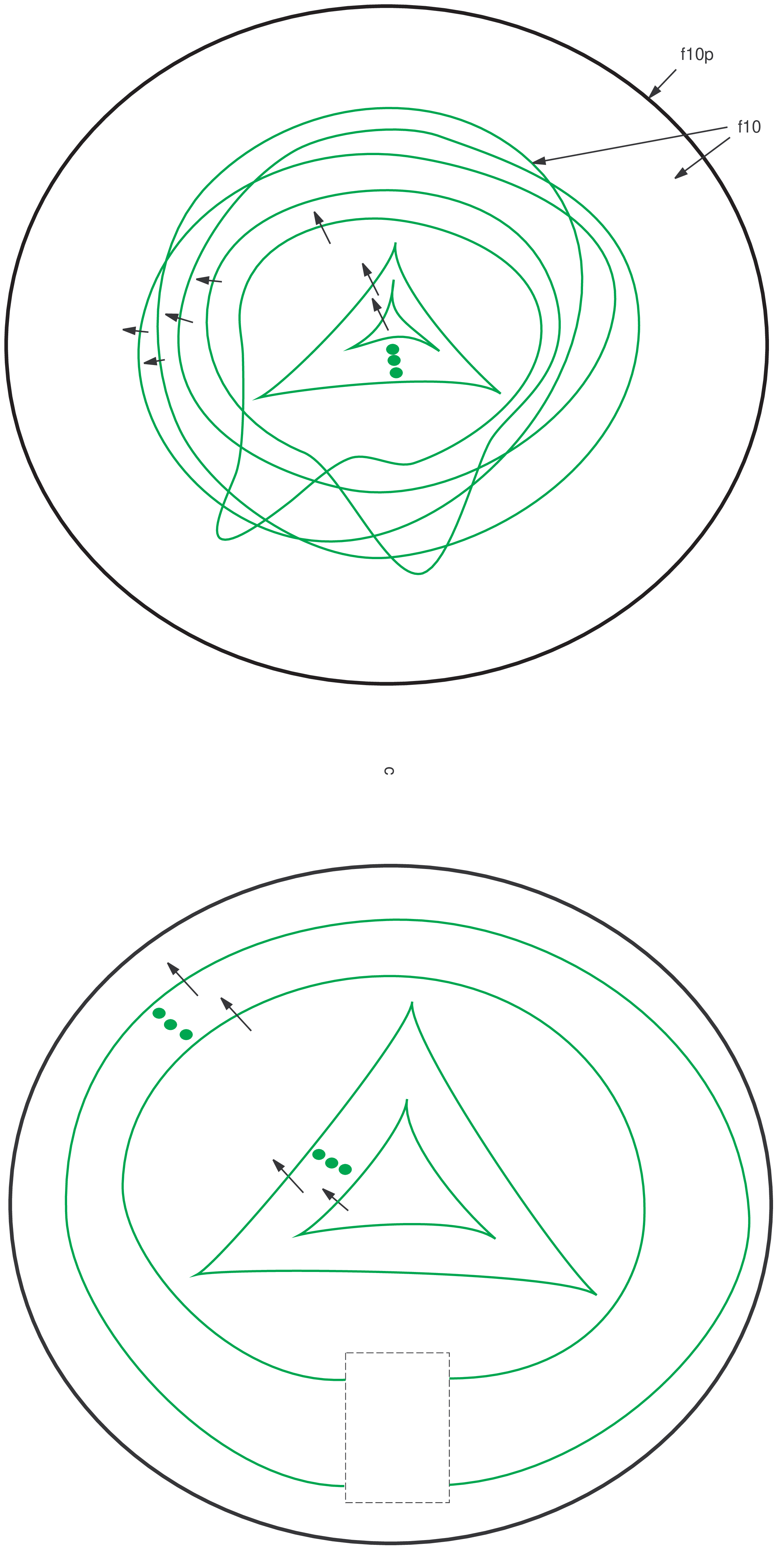}
\caption{Making a family of concentric ``triangles''. The
square box in the bottom represents a Cerf graphic without cusps
in the sense of \cite{GK4}.}
\label{fig8}
\end{figure}

Then, by the same procedure as in the last
part of the proof of \cite[Theorem~7.1]{BS1},
we turn a Lefschetz singularity image to
a ``triangle'' with three indefinite cusp images by the wrinkle 
move, push the rest of the Lefschetz singularity images
into the interior of the ``triangle'', and repeat 
the process. Then we get a generic map
$f_{10}: M \to \Delta_1 \cong D^2$,
which is easily seen to be
a desired relative trisection map (see Fig.~\ref{fig8}).
This completes the proof.
\end{proof}

\begin{rk}
In the case where $M$ is a closed connected orientable $4$--dimensional manifold,
a \emph{simplified} trisection map is known to exist on $M$
in the sense that the Cerf graphic is assumed to be trivial without
any cusps nor crossings
(see \cite{BS1}). In the case where $\partial M \neq
\emptyset$, we do not know if such a simplified
relative trisection map always exists or not.
\end{rk}

\begin{rk}
Even if for a component of
$\partial M$, the open book structure induced by
$OB$ has empty binding, the above proof works.
\end{rk}

The following is due to \cite{GK4}, see also \cite{CGP0, CGP}.

\begin{cor}\label{cor1}
Every compact connected orientable $4$--dimensional
manifold $M$ with nonempty boundary $\partial M$, which
may be disconnected, admits a relative trisection.
Furthermore, given any open book structure $OB$ on
$\partial M$, such a relative trisection can be chosen
so that it is consistent with $OB$.
\end{cor}

\begin{figure}[t]
\centering
\psfrag{D}{$D^2$}
\includegraphics[width=0.95\linewidth,height=0.3\textheight,
keepaspectratio]{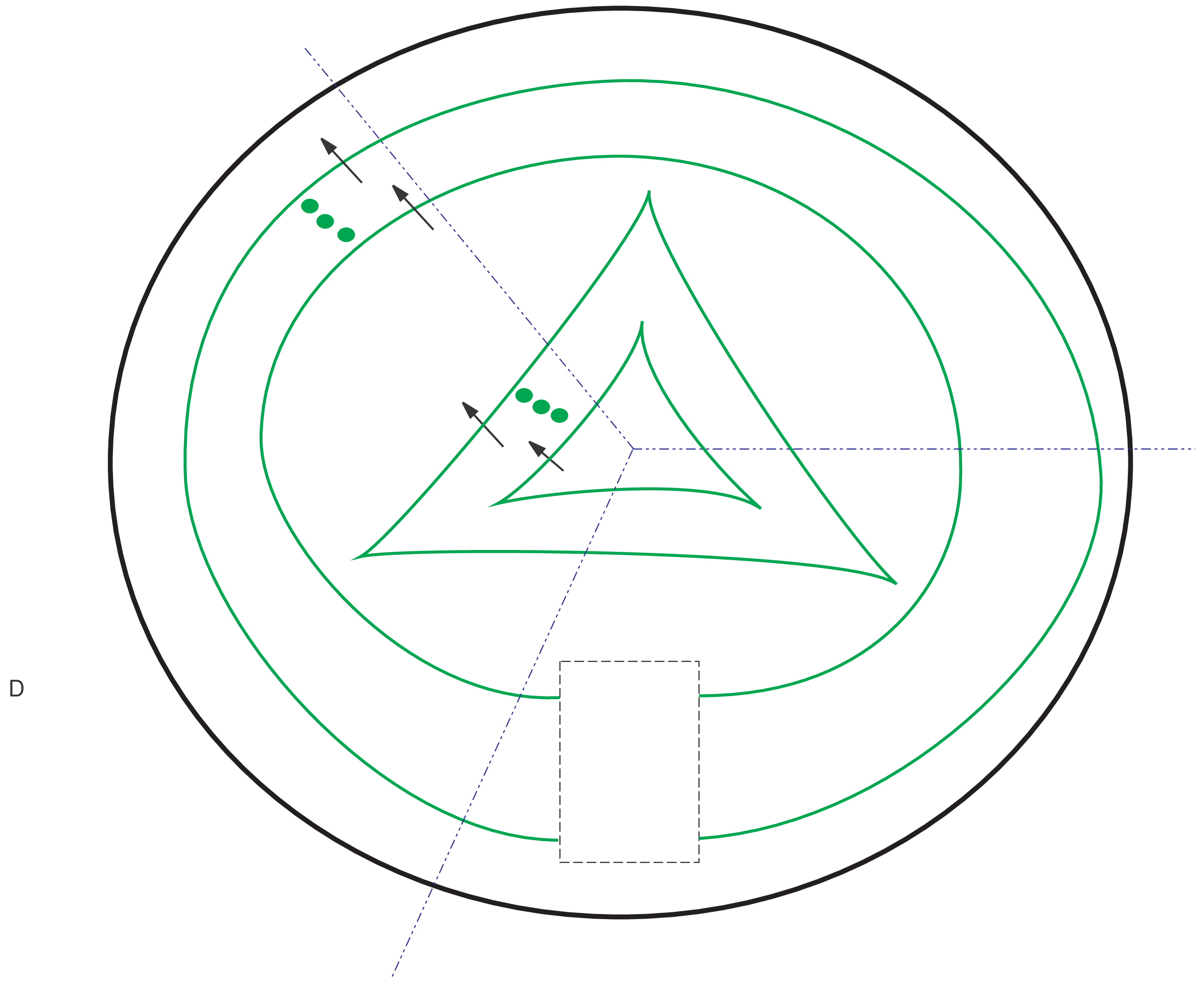}
\caption{Three sectors in $D^2$ whose inverse images give
a relative trisection of $M$.}
\label{fig14}
\end{figure}

For the definition of a relative trisection, the reader
is referred to \cite{CGP0, CGP, GK4}. Such a trisection
can be obtained by taking the inverse images
of the three sectors as depicted in Fig.~\ref{fig14}.
In fact, in the literature, relative trisections have been formulated
for the case where the binding on each component of $\partial M$
is nonempty \cite{CGP0, CGP}.
On the other hand, in \cite{GK4}, there is no
such assumption, although the boundary of the $4$--manifold
is assumed to be connected.
In the most general case, we can formulate a relative
trisection in a similar way: however, each piece
may not be diffeomorphic to a $4$--dimensional
handlebody.

Note that for a closed connected orientable $4$--dimensional manifold $M$,
Gay--Kirby \cite{GK4} showed the existence of a 
\emph{trisection map} on $M$ and as a result
proved the existence of a trisection on $M$.
(They also gave another proof based on handlebody
decomposition of $M$.)
Our result is along the same line of an idea, but
the proof is quite different from the one given in \cite{GK4}.

\section{Realizing given relative trisection by a map}\label{section3.5}

In this section, we give another proof to Theorem~\ref{thm1}
based on the fact that every compact connected orientable
$4$--manifold with an open book structure $OB$ given on the boundary
admits a relative trisection consistent with $OB$.
More precisely, we show the following.

\begin{thm}\label{thm2}
Let $M$ be a compact connected orientable $4$--dimensional
manifold with nonempty boundary, which may be disconnected, such that
an open book structure $OB$ is given on $\partial M$.
We assume that the binding of the open book structure induced
by $OB$ on each component of $\partial M$ is nonempty.
Suppose that $M$ admits a relative trisection $(X_1, X_2, X_3)$
consistent with $OB$. Then, there exists a relative
trisection map $f : M \to D^2$ compatible with $OB$ which
induces the relative trisection $(X_1, X_2, X_3)$.
\end{thm}

\begin{rk}
For closed $4$--manifolds, a corresponding result
seems to be folklore. For example, the fact
is stated in \cite[\S2.2]{Hayano1} without proof.
\end{rk}

In the following, we use the following generalization
of a $3$--dimensional handlebody (see, for example, \cite{Goda}).

\begin{defn}
Let $P$ be a compact orientable possibly disconnected surface with
possibly nonempty or disconnected boundary.
Let $C$ be a compact connected orientable
$3$--dimensional manifold obtained from $P \times [0, 1]$ by attaching
a finite number of $3$--dimensional $1$--handles along $P \times \{1\}$.
We call such a $3$--manifold $C$ a \emph{sutured compression body},
where $\partial P \times [0, 1]$ is the \emph{suture}.
Furthermore,
we set $\partial_- C = P \times \{0\}$
and we denote by $\partial_+ C$ the closure of $\partial C
\setminus (\partial_- C \cup (\partial P \times [0, 1]))$.

When $P = \emptyset$, we start with a $3$--dimensional disk $D^3$ and attach
$1$--handles along its boundary. Then, $C$ is nothing but a usual
$3$--dimensional (orientable) handlebody.
\end{defn}

\begin{proof}[Proof of Theorem~\textup{\ref{thm2}}]
In the following, we use terminologies and notations in \cite{CGP}.
Suppose that the given decomposition 
$X = X_1 \cup X_2 \cup X_3$ is a 
$(g, k; p_1, b_1; p_2, b_2: \ldots ; p_m, b_m)$--trisection, where $g$ is the
genus of the surface $\Sigma = X_1 \cap X_2 \cap X_3$,
$X_i$ are diffeomorphic to $\natural^k (S^1 \times D^3)$,
$i = 1, 2, 3$, and the page of the open book structure
on the $j$--th component of $\partial M$ has genus $p_j$ and
$b_j \geq 1$ boundary components, $j = 1, 2, \ldots, m$. 

Let $(\Sigma, \alpha, \beta, \gamma)$ be
a relative trisection diagram associated with the
relative trisection $(X_1, X_2, X_3)$
in the sense of \cite{CGP}.
Let $C_\alpha$ be the sutured compression body
obtained from $\Sigma \times [0, 1]$ by attaching
$3$--dimensional $2$--handles along $\alpha \subset \Sigma
\times \{1\}$. We construct $C_\beta$ and $C_\gamma$ similarly.
Let $\Delta$ be a triangle with edges $e_\alpha$,
$e_\beta$ and $e_\gamma$ as in \cite[Proof of Theorem~3]{CGP}.
Then, we can construct a smooth map
$f_1 : M_1 \to \R^2$ whose image is as depicted in
Fig.~\ref{fig16}, where $M_1$ is the $4$--manifold
obtained from $\Delta \times \Sigma$ by attaching
$e_\alpha \times C_\alpha$,
$e_\beta \times C_\beta$ and
$e_\gamma \times C_\gamma$ along
$e_\alpha \times (\Sigma \times \{0\})$, 
$e_\beta \times (\Sigma \times \{0\})$
and $e_\gamma \times (\Sigma \times \{0\})$, respectively,
$\ell = \sum_j (2p_j+b_j-1)$, 
$n = k - \ell$, $s = g-n-\sum_j p_j$,
$f_1$ has folds as its singularities,
and the green curves depict its singular value set.
Such a smooth map $f_1$ can be constructed by using
the projection $\Delta \times \Sigma \to \Delta$
and Morse functions $C_\alpha \to [0, 1]$,
$C_\beta \to [0, 1]$ and $C_\gamma \to [0, 1]$
consistent with the attached $2$--handles as described above.
Note that the final $m-1$ $2$--handles are attached along separating
simple closed curves for each of $C_\alpha, C_\beta$ and $C_\gamma$.
The inverse images of the red parts in Fig.~\ref{fig16}
correspond to
$M_{\alpha \beta}$, $M_{\beta \gamma}$, and $M_{\alpha \gamma}$
(in $\partial M_1$)
in the notation in \cite[Proof of Theorem~3]{CGP}.

\begin{figure}[t]
\centering
\psfrag{ea}{$e_\alpha$}
\psfrag{eb}{$e_\beta$}
\psfrag{ec}{$e_\gamma$}
\psfrag{D}{$\Delta$}
\psfrag{ca}{$e_\alpha \times C_\alpha$}
\psfrag{cb}{$e_\beta \times C_\beta$}
\psfrag{cc}{$e_\gamma \times C_\gamma$}
\psfrag{Mab}{$M_{\alpha \beta}$}
\psfrag{Mbc}{$M_{\beta \gamma}$}
\psfrag{Mca}{$M_{\gamma \alpha}$}
\psfrag{s}{$s$}
\psfrag{l}{$n+m-1$}
\includegraphics[width=0.95\linewidth,height=0.4\textheight,
keepaspectratio]{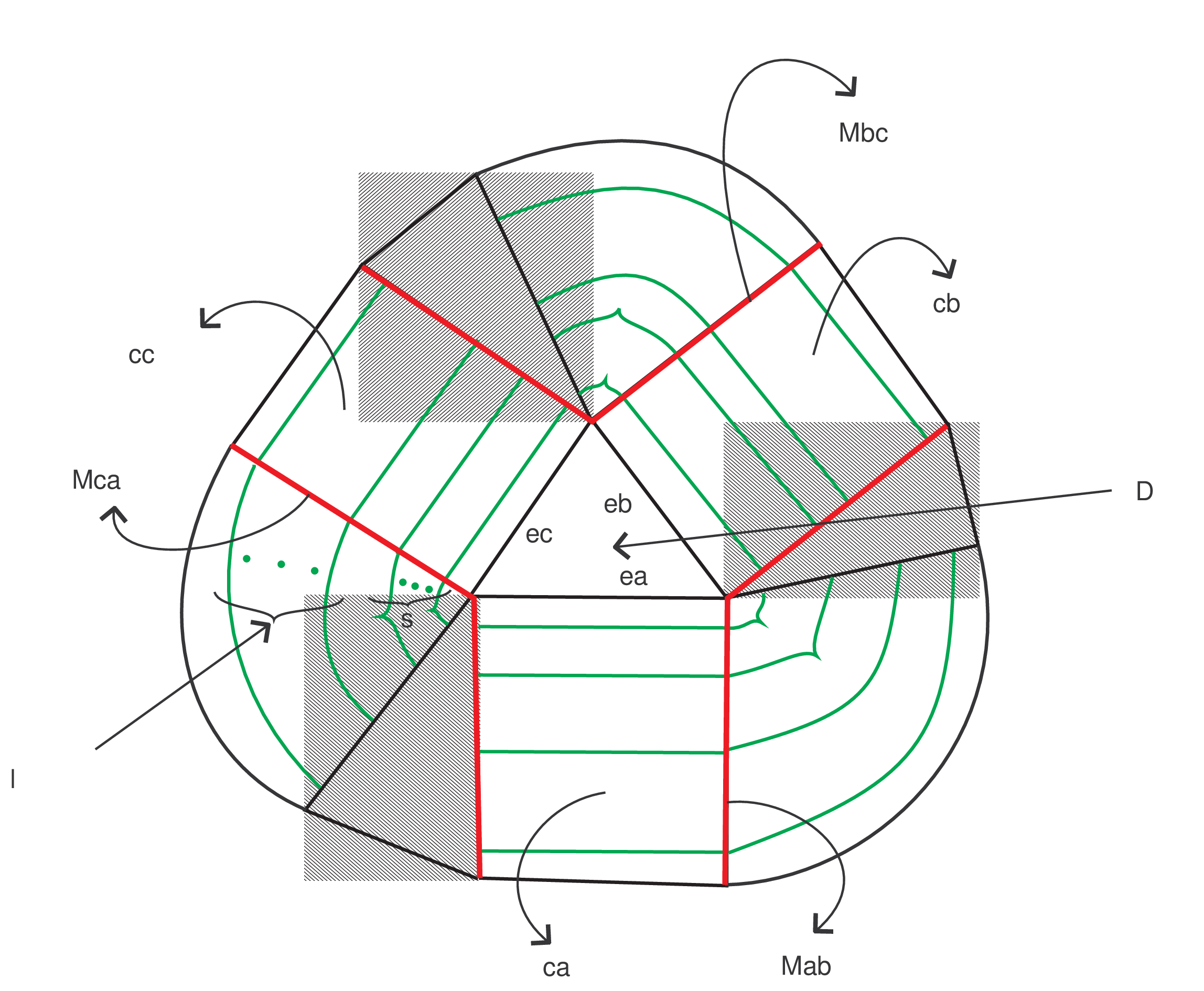}
\caption{Image of a generic map}
\label{fig16}
\end{figure}

Now, since $(\Sigma, \alpha, \beta)$ is handleslide
diffeomorphic to the standard one, we can find
a generic smooth map into a shaded triangular region 
adjacent to the intersection of the images of 
$e_\beta \times C_\beta$ and
$M_{\alpha \beta}$, extending $f_1$, by
using an appropriate $1$--parameter family
of Morse functions on $C_\beta$ corresponding to handleslides, 
in such a way that it represents a Cerf graphic without cusps,
and the resulting map can be naturally further extended
to a generic map into $\R^2$ so that
its singular values are as depicted in Fig.~\ref{fig16}.
Similar arguments apply to $(\Sigma, \beta, \gamma)$
and $(\Sigma, \gamma, \alpha)$.
Then, according to \cite{CGP0, CGP}, the source
$4$--manifold is diffeomorphic to $M = X_1 \cup
X_2 \cup X_3$ and the resulting generic map
is the desired relative trisection map.
This completes the proof.
\end{proof}

When the bindings of the open book
structures on some components of $\partial M$ are empty,
each $X_i$ is not diffeomorphic to a $4$--dimensional
handlebody, but is obtained from $P \times D^2$ by
attaching some $4$--dimensional $1$--handles along
$P \times \partial D^2$, where $P$ is a page of $OB$ and
may have closed components.
In such a case, we do not know if the conclusion of
Theorem~\ref{thm2} holds or not.
See also the paragraph just after Corollary~\ref{cor1}.

\section{A decomposition of $4$--manifold}\label{section4}

In this section, we show that a compact connected oriented $4$--manifold
always decomposes as the union of a PALF and a sutured compression body
bundle over $S^1$ compatible with a given open book structure
on the boundary.

Let $M$ be a compact connected oriented $4$--dimensional
manifold such that an open book structure $OB$ is
given on $\partial M \neq \emptyset$.

\begin{defn}
Let $f : M \to D^2$ be a smooth map compatible
with an open book structure $OB$ on $\partial M$.
A positive or negative Lefschetz singularity of $f$
is \emph{allowable} if its vanishing cycle is
homologically non-trivial in a regular fiber of $f$.
We say that $f$ is a \emph{positive allowable
Lefschetz fibration} (\emph{PALF}, for short)
if $f$ is a LF such that the fiber 
has nonempty boundary and its
Lefschetz singularities are all allowable.

Note that if $f$ is a PALF, then its restriction
to $f^{-1}(\partial D^2)$ is a surface bundle over the
circle and the monodromy is a composition of some
positive Dehn twists.
Furthermore, both $\partial M (\neq \emptyset)$ and the fibers are 
connected.
\end{defn}

It is known that a compact connected oriented
$4$--dimensional manifold $M$ with nonempty
boundary is a Stein surface (with strictly pseudoconvex boundary)
if and only if
it admits a PALF $f : M \to D^2$ \cite{AO, LP}.
(For the definition of a Stein surface, see \cite{LP}.)

\begin{lem}\label{lem1}
The Lefschetz singularities that appear in the
proof of Theorem~\ref{thm1} are all allowable,
in the sense that their vanishing cycles
are homologically non-trivial on the regular fibers.
\end{lem}

\begin{proof}
We observe that $f_4 : M \to S^2$ has only
indefinite folds and cusps and has no Lefschetz
singularities. On the other hand,
$f_5$ may have Lefschetz singularities. They all
appear in the process of applying \cite[Theorem~4.1]{BS1},
in which only the unsink and push moves involve
Lefschetz singularities. By \cite[Figure~1]{Lek},
we see that the Lefschetz singularity that appears
in the process of an unsink move is allowable,
since for a vanishing cycle $\gamma$ on a regular fiber, 
there is another simple closed curve $\gamma'$ on the
same fiber that intersects $\gamma$ transversely
at one point.
Furthermore, such a pair of curves persist even after
a push move.
Therefore, the Lefschetz singularities of $f_5$
are all allowable.

In the transition from $f_5$ to $f_6$, some new
Lefschetz singularities my appear and they
are all allowable by the same reason.
For $f_7, f_8$ and $f_9$, their Lefschetz singularities
are essentially the same as those of $f_6$ and
are allowable. This completes the proof.
\end{proof}

We have the following corollary to Theorem~\ref{thm1}.

\begin{cor}\label{cor2}
Let $M$ be a compact connected oriented $4$--dimensional
manifold with nonempty boundary $\partial M$, which
may be disconnected.
Suppose an open book structure $OB = (\partial M, B, T(B), \varphi)$ 
on $\partial M$ with $B \neq \emptyset$ is given.
Then, the fibration $\partial M \setminus \Int{T(B)} \to S^1$
induced by $\varphi$
can be extended to a smooth fiber
bundle $M_C \to S^1$ with fiber a sutured compression body $C$
with $\partial_- C$ being
a page of $OB$ and with
$\partial_+ C$ being
a compact connected oriented surface $\Sigma$
such that 
the closure $M_{\mathrm{PALF}}$ of $M \setminus M_C$
admits a PALF over a closed $2$--disk $D^2$ (with
fiber being compact connected oriented
surface $\Sigma$ with boundary) consistent
with the $\partial_+ C$--bundle structure
induced by the fibration $M_C \to S^1$ over $\partial D^2$.
\end{cor}

\begin{proof}
Let us consider the map $f_9 : M \to D^2$
constructed in the proof of Theorem~\ref{thm1}
(refer to Fig.~\ref{fig7}).
First, if there are no Lefschetz singularities, then
we can apply the birth move as in the top-left
of \cite[Fig.~3]{BS1}, and then the unsink moves
twice so that we have two (positive or negative)
Lefschetz singularities. Hence, we may assume
that there are at least one Lefschetz singularities.

By \cite[\S6]{Lek}, one can
replace each negative Lefschetz singularity
by three positive Lefschetz singularities together
with a circle of indefinite fold whose image is
outward directed. By applying this procedure
for each negative Lefschetz singularity and
by using $\mathrm{R2}^0$ moves together with
push moves as in \cite{BS1}, we may assume that
the map on $M_{\mathrm{PALF}} = (f_9)^{-1}(\Delta_3)$ 
is a Lefschetz fibration with positive
Lefschetz singularities. Note that they
are all allowable by Lemma~\ref{lem1}
and that the fibers are connected with nonempty boundary.

Set $M_C = f^{-1}(D^2 \setminus \Int{\Delta_3})$.
Then, $f_9|_{M_C} : M_C \to D^2 \setminus \Int{\Delta_3}$
composed with the radial projection 
$D^2 \setminus \Int{\Delta_3} \to \partial D^2$
is a desired fibration.
This completes the proof.
\end{proof}

\begin{cor}\label{cor3}
Let $M$ be a compact connected oriented $4$--dimensional
manifold with nonempty boundary $\partial M$, which
may be disconnected.
Suppose an open book structure $OB$ on $\partial M$
is given. Then, there always exists a BLF
$M \to D^2$ compatible with $OB$.
\end{cor}

Corollary~\ref{cor3} follows immediately from the above proof
(without assuming that the binding of $OB$ is nonempty).

We do not know if $M$ admits a \emph{near-symplectic structure}
compatible with the BLF in a certain sense.
For near-symplectic structures, the reader is
referred to \cite{ADK, BS1}.

Note that the binding $B$ of $OB$ is contained in $M_{\mathrm{PALF}}$
and that $M_C$ is a neighborhood of $\partial M \setminus \Int{T(B)}$
in Corollary~\ref{cor2}.
The above corollaries mean, more or less, that
the open book structure $OB$ given on $\partial M$
can be extended to $M$ except for some
positive (allowable) Lefschetz singularities.
The fiber of the Lefschetz
fibration over the disk is connected, since
the $3$--dimensional fiber $C$ of the fibration $M_C \to S^1$
is connected and is constructed from the product of
a page and an interval by attaching some $3$--dimensional
$1$--handles.
This argument also shows that $C$ is a sutured compression body.

Note that in general,
the total space of a PALF
naturally carries a symplectic structure (in fact,
a Stein structure) which has strictly pseudoconvex boundary. Therefore,
$M_{\mathrm{PALF}}$ admits such a symplectic structure.

Corollary~\ref{cor3} strengthens \cite[Corollary~5.3]{GK1}
in several ways. First, the $4$--dimensional
manifold $M$ may not necessarily have the handlebody
decomposition by $0$--, $1$-- and $2$--handles.
Second, the open book structure on $\partial M$ may not
support an overtwisted contact structure.
Third, we do not need to stabilize the open book
structure $OB$ on the boundary for its extension.

\begin{rk}
Note that in Corollary~\ref{cor2},
$M$ admits a handlebody decomposition
such that $M_{\mathrm{PALF}}$ has a
handlebody decomposition consisting of one $0$--handle,
some $1$-- and $2$--handles, and that $M$ is obtained by
attaching some $2$-- and $3$--handles to $M_{\mathrm{PALF}}$.
This follows from the construction together with
the fact that the indefinite fold image of the BLF
is outward directed.
\end{rk}

Note also that the boundary of $M_{\mathrm{PALF}}$ is Stein fillable.

We also have the following as a byproduct.

\begin{cor}\label{cor:right}
Let $P$ be a compact (possibly disconnected) oriented surface
(possibly with boundary)
and $h : P \to P$ an orientation preserving diffeomorphism which
fixes $\partial P$ pointwise. Then, there exists a $3$--dimensional sutured
compression body $C$ with $\partial_- C = P$ and with $\partial_+ C$ connected
together with an orientation
preserving diffeomorphism $H : C \to C$ which fixes the suture
pointwise such that $H|_{\partial_+ C}$ is isotopic to
a composition of finitely many positive Dehn twists.
\end{cor}

\begin{proof}
Using $P$ as page and $h$ as monodromy, we can construct a (possibly disconnected)
closed orientable $3$--dimensional manifold $N$ with an open book structure $OB$
whose monodromy coincides with $h$.
By \cite{R}, there exists a compact orientable $4$--dimensional manifold $M$ with
$\partial M = N$. By taking interior connected sum, we may assume that
$M$ is connected.
Then, by Corollary~\ref{cor2}, we have a sutured compression
body $C$ and a smooth $C$--bundle $M_C$ over $S^1$ whose monodromy
$H$ extends $h$. Then $H|_{\partial_+ C}$ is isotopic to a composition
of finitely many positive Dehn twists, since it 
coincides with the monodromy of the surface bundle over
the boundary of a $2$--disk of a PALF.
This completes the proof.
\end{proof}

Note that not every orientation preserving
diffeomorphism of a compact oriented surface is isotopic to a composition
of finitely many positive Dehn twists (for example, see \cite{HKM}).

\section*{Acknowledgment}\label{ack}
The author would like thank Toru Kojima for useful
discussions about relative trisections.
He would also like to thank Professor \.{I}nan\c{c} Baykur
for stimulating discussions which motivated the
theme of this paper. 
The author would also like to thank Professor Vincent
Blanl\oe il and the Advanced Mathematical Research Institute
of the University of Strasbourg
for the hospitality
during the preparation of the manuscript.
This work has been supported in part by JSPS KAKENHI Grant Numbers 
JP22K18267, JP23H05437, JP25K24558. 

\section*{Declaration on the use of AI tools}

Generative AI tools (specifically Microsoft Copilot and Google Gemini 3.8 Flash) 
were used during the preparation of this manuscript solely for the purpose of literature search. 
The author reviewed and verified all sourced references, and accept full 
responsibility for the accuracy.

\end{document}